\documentclass[letterpaper, 10 pt, conference]{ieeeconf}
\IEEEoverridecommandlockouts
\usepackage{graphics}
\usepackage{subfigure}
\usepackage{epsfig}
\usepackage{mathptmx}
\usepackage{times}
\usepackage{amsfonts}
\usepackage{amsmath}
\usepackage{amssymb}
\usepackage{mathtools}
\usepackage{mathrsfs}
\usepackage{bm}
\usepackage{algorithm}
\usepackage{commath}
\usepackage{algpseudocode}
\usepackage{romannum}
\newtheorem{theorem}{Theorem}
\DeclareMathAlphabet{\mathcal}{OMS}{cmsy}{m}{n}
\usepackage{multirow}
\usepackage[colorlinks,
        linkcolor=black,
        anchorcolor=black,
        citecolor=black,
        urlcolor=black
        ]{hyperref}
\usepackage{ifthen}
    \usepackage[noprefix]{nomencl}
\usepackage{cite}
\usepackage{multirow}
\usepackage{multicol}
\usepackage{booktabs}
\usepackage{footnote}
\usepackage{arydshln}

\title{{\LARGE \bf
    Distributionally Robust Resource Allocation with Trust-aided Parametric Information Fusion
}}

\author{
    Yanru Guo, Bo Zhou, Ruiwei Jiang, Xi (Jessie) Yang, Siqian Shen
    \thanks{
   The authors are with Department of Industrial and Operations Engineering, University of Michigan, Ann Arbor, MI 48109, US.
        (Emails: \{yanruguo, bozum, ruiwei, xijyang, siqian\}@umich.edu)
    }
}

\begin{document}

    \maketitle
    \thispagestyle{empty}
    \pagestyle{empty}

\begin{abstract}
    Reference information plays an essential role for making decisions under uncertainty, yet may vary across multiple data sources. In this paper, we study resource allocation in stochastic dynamic environments, where we perform information fusion based on trust of different data sources, to design an ambiguity set for attaining distributionally robust resource allocation solutions. We dynamically update the trust parameter to simulate the decision maker's trust change based on losses caused by mis-specified reference information. We show an equivalent tractable linear programming reformulation of the distributionally robust optimization model and demonstrate the performance in a wildfire suppression application, where we use drone and satellite data to estimate the needs of resources in different regions. We demonstrate how our methods can improve trust and decision accuracy. The computational time grows linearly in the number of data sources and problem sizes.
    
\end{abstract}

\section{Introduction}\label{1}

    In the realms of stochastic optimization and system control, decision-making challenges involve information uncertainty, which often stems from limited data and imprecise measurements. Over the last decade, distributionally robust optimization (DRO) has been widely used for attaining the best worst-case performance under ambiguously known uncertain parameters \cite{delage2010distributionally, jiang2016data, mohajerin2018data}. In broader resource allocation problems, \cite{yang2021multi} formulates a multi-period DRO model to dynamically optimize pre-positioning of emergency supplies under demand uncertainty; \cite{wang2021two} uses DRO to optimize resource allocation in disaster relief; \cite{basciftci2023resource} compares a stochastic programming (SP) approach with DRO for epidemic control resource allocation under stochastic spatiotemporal demand, and shows that the DRO approach can limit the number of unvaccinated or untested people by paying higher cost. 
   
    Building on the emergence of DRO models, in this paper, we consider stochastic resource allocation with unknown demand whose distribution can be inferred from multiple data sources via data fusion \cite{durrant2016multisensor}. 
      The existing techniques for fusing data can be classified into two main categories: probability-based methods \cite{abdulhafiz2014bayesian, dietrich2017probabilistic, noyvirt2012human} and Artificial Intelligence (AI)-based methods \cite{howcroft2017prospective, pan2020human, bloomfield2021machine, nsengiyumva2022critical}, but not all are adaptable for characterizing the ambiguity set in DRO models. In our approach, we use trust \cite{lee2004trust} as weights to fuse information from different sources, and its dynamic nature \cite{guo2021modeling, yang2023toward} allows us to update trust on different sources over time.

     Our work complements the existing studies in stochastic resource allocation using  DRO with trust being the ambiguously known uncertain parameter and we take multi-source reference information into consideration. 
    
    The main contributions of the paper are three-fold: 
    \begin{itemize}
        \item We develop a multi-reference distributionally robust optimization (MR-DRO) model for stochastic resource allocation. We combine predicted information from multiple sources to construct the ambiguity set.
        \item We design a trust update process to simulate trust variation over time once more data becomes available. We utilize historical data and outcomes to select proper trust that adapts to the relative prediction errors. 
        \item We show, via extensive computational results, that the MR-DRO model yields better results than the ones obtained by fully trusting a single source.
    \end{itemize}

    The remainder of the paper is organized as follows. In Section \ref{2}, we formulate the MR-DRO model by constructing the trust-aided ambiguity set and derive its tractable reformulation. In Section \ref{3}, we propose a trust update process to simulate trust change over time. Section \ref{4} includes detailed numerical experiments and result analysis. In Section \ref{5}, we conclude the paper and state future research directions.

\section{Models and Solution Approaches}\label{2}
    
    We use a wildfire suppression example to explain our models and algorithms, where a decision maker predicts suppression resource demand from multiple data sources, such as satellites and drones, for a geographic region. We consider an ambiguously known distribution of the demand and hedge against the risk of inappropriate allocation that may result in over-served or unmet demand in each region. We further assume that the decision maker has different trust on different sources, depending on their historical feedback and accuracy. 
        
    \subsection{Problem Description and Notation}\label{2.1}
        Let $K$ be the number of regions. Let $\bm{c}^{u} = [c^{u}_{1},\ldots,c^{u}_{K}]^\mathsf{T}$ and $\bm{c}^{o} = [{c}^{o}_{1},\ldots,{c}^{o}_{K}]^\mathsf{T}$ be unit penalty costs of unmet and over-served demand in each region, respectively.  (Throughout the paper, we use bold symbols to denote the vector form of a decision variable or a parameter.) We have an overall resource budget $B > 0$. The demands in all regions are uncertain and are denoted by a random vector \(\bm{\xi} = [\xi_{1},\ldots,\xi_{K}]^\mathsf{T}\), where $\bm{\xi} \in \mathbb{R}^{K}_{+}$ obeys a probability distribution $\mathbb{P}$. We define decision variable $x_{k} \geq 0$ for all $k \in [K]$ as the amount of resource allocated to region $k$, where $[K] = \{1,\ldots,K\}$. 

        \subsection{SP and DRO Models}\label{2.2}
        Considering fully known distribution $\mathbb P$ of uncertain demand $\xi$, one can minimize the expected total cost of over-served and unmet demand in all regions via an SP model, i.e., 
            \begin{subequations}\label{eq:SO}
                \begin{gather}
                    \inf_{\bm{x} \in \mathbb{X}}\mathop{\mathbb{E}}\limits_{\bm{\xi}\sim\mathbb{P}}\left[(\bm{c}^{u})^{\mathsf{T}}(\bm{\xi}-\bm{x})^{+}+(\bm{c}^{o})^{\mathsf{T}}(\bm{x}-\bm{\xi})^{+}\right]\label{eq:SO-obj}\\
                    \textrm{s.t.}~\mathbb{X} = \left\{\bm{x} \in \mathbb{R}^{K}_{+}: \sum^{K}_{k=1}x_{k} \leq B\right\}\label{eq:SO-budget}.
                \end{gather}
            \end{subequations}
      
            To solve \eqref{eq:SO}, one can apply the Monte Carlo sampling approach to replace $\mathbb{P}$ with an empirical distribution constructed by $|\Omega|$ scenarios, each having an equal probability $p^{\omega} = 1/|\Omega|$. The resulting problem is called the Sample Average Approximation (SAA) problem \cite{kleywegt2002sample}. 
            Specifically, for each scenario $\omega \in \Omega$, we denote $\bm{\xi}^{\omega} = (\xi^{\omega}_{1},\ldots,\xi^{\omega}_{K})^{\mathsf{T}}, \bm{\xi}^{\omega} \in \mathbb{R}^{K}$ as the demand realization in scenario $\omega$ and reformulate \eqref{eq:SO} as
            \begin{equation}\label{eq:SP}
                \min_{x \in \mathbb{X}}\quad\sum_{\omega \in \Omega}p^{\omega}\left((\bm{c}^{u})^{\mathsf{T}}(\bm{\xi}^{\omega}-\bm{x})^{+}+(\bm{c}^{o})^{\mathsf{T}}(\bm{x}-\bm{\xi}^{\omega})^{+}\right).
            \end{equation}
            
            However, in practice, the demand distribution $\mathbb{P}$ is hard to be known precisely and the acquisition of samples might be impossible or expensive. 
            In light of this issue, we derive a variant of \eqref{eq:SO} as a DRO model that  accounts for the ambiguity of the true distribution $\mathbb{P}$. For notation simplicity,  we let $\ell(\bm{x},\bm{\xi}) = (\bm{c}^{u})^{\mathsf{T}}(\bm{\xi}-\bm{x})^{+}+(\bm{c}^{o})^{\mathsf{T}}(\bm{x}-\bm{\xi})^{+}$ and consider:
            \begin{equation}\label{eq:DRO}
                    \inf_{\bm{x} \in \mathbb{X}}\left\{\sup_{\mathbb{P} \in \mathcal{P}}\mathop{\mathbb{E}}\limits_{\bm{\xi}\sim\mathbb{P}}[\ell(\bm{x},\bm{\xi})]\right\},
            \end{equation}
            where the definition of ambiguity set $\mathcal{P}$ is given by: 
            \begin{equation}\label{eq:DRO-set}
                \mathcal{P}:=\left\{\mathbb{P} \in \mathcal{M}(\Xi) : d_{W}(\mathbb{P},\hat{\mathbb{P}}_{N}) \leq \epsilon\right\}.
            \end{equation}
The set $\mathcal{P}$ is a Wasserstein ambiguity set  \cite{mohajerin2018data} centered at the empirical distribution $\hat{\mathbb{P}}_{N}$ and bounded by radius $\epsilon$, constructed based on the Wasserstein metric using $L_1$-norm. We denote $\Xi$ as the Cartesian product of closed convex sets $\Xi_{k}$, where $\xi_{k} \in \Xi_{k}$ for all $k \in [K]$. We define the Wasserstein metric on the space $\mathcal{M}(\Xi)$ of all probability distributions $\mathbb{P}$ supported on $\Xi$ with $\mathbb{E}^{\mathbb{P}}[\|\bm{\xi}\|] = \int_{\Xi} \|\bm{\xi}\|\mathbb{P}(d\bm{\xi}) < \infty$ and $d_{W}: \mathcal{M}(\Xi) \times \mathcal{M}(\Xi) \rightarrow \mathbb{R}_{+}$ via:
            \begin{equation}\label{eq:DRO-WassersteinMetric}
                d_{W}(\mathbb{P},\hat{\mathbb{P}}_{N}) = \inf \left\{ \int_{\Xi^{2}}\|\bm{\xi}-\hat{\bm{\xi}}\|\Pi(d\bm{\xi},d\hat{\bm{\xi}})\right\},
            \end{equation}
            where $\|\cdot\|$ represents $L_1$-norm, and $\Pi$ is the joint distribution of $\bm{\xi}$ and $\hat{\bm{\xi}}$ with marginal distributions $\mathbb{P}$ and $\hat{\mathbb{P}}_{N}$, respectively and $\mathbb{P}, \hat{\mathbb{P}}_{N} \in \mathcal{M}(\Xi)$.
            
            The objective of DRO model \eqref{eq:DRO} is to minimize the maximum expected costs/losses over all possible distributions in the ambiguity set $\mathcal{P}$. Therefore, we will carefully select the radius \(\epsilon\) so that $\mathcal{P}$ can cover the true distribution \(\mathbb{P}\) with sufficiently high probability, and meanwhile, it will not result in overly conservative solutions. 

        \subsection{Trust-Aided Parametric Data-fusion Ambiguity Set}
        \label{2.3}
        
\begin{figure}[!htbp]
                \centering
                \includegraphics[width=0.35\textwidth]{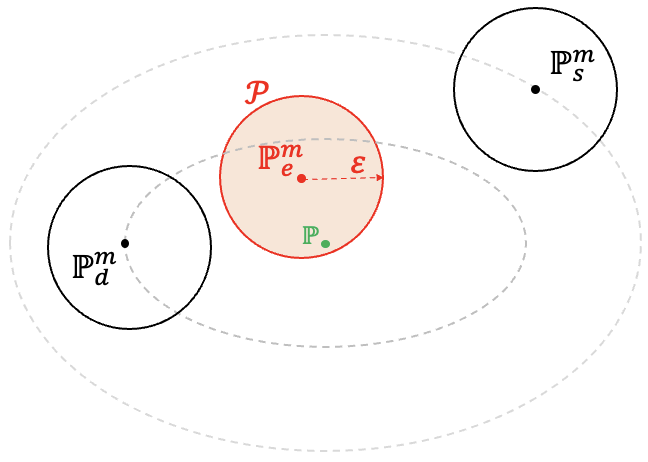}
                \caption{Illustration of a parametric data-fusion trust-aided ambiguity set ($\mathbb{P}$: true demand distribution; $\mathbb{P}^{m}_{s}$: predicted distribution provided by the satellite; $\mathbb{P}^{m}_{d}$: predicted distribution provided by the drone; $\mathbb{P}^{m}_{e}$: empirical distribution; $\mathcal{P}$: ambiguity set)}
                \label{fig:parametric-fusion}
            \end{figure}

        Next, we discuss how to obtain the empirical distribution $\hat{\mathbb{P}}_{N}$. Denote $M$ as the number of wildfire events, and at the beginning of each event $m \in [M]$, we receive the predicted distribution of demands/wildfire status from both satellite and drone, denoted as $\mathbb{P}^{m}_{s}$ and $\mathbb{P}^{m}_{d}$, respectively. We consider the trust held for drone data after event $m$ as $\bm{t}^{m} = (t^{m}_{1},\ldots,t^{m}_{K})^{\mathsf{T}}$ and thus the trust on satellite data is $(1-t^m_k)$ for each region $k \in [K]$. We conduct information fusion in a parametric way as follows. Assume that $\hat{\mathbb{P}}_{N}$ in \eqref{eq:DRO-set} consists of $N$ i.i.d. samples from a Normal distribution $\mathbb{P}^{m}_{e} \approx \mathcal{N}(\bm{\mu}^{m},\Sigma^{m})$. For $\bm{\mu}^{m} = [\mu^{m}_{1},\ldots,\mu^{m}_{K}]^\mathsf{T}$, let $\mu^{m}_{k} = t^{(m-1)}_{k}\mu^{m}_{dk} + (1-t^{(m-1)}_{k})\mu^{m}_{sk}$, $\forall k \in [K]$, and for $\Sigma^{m} = \text{diag}((\sigma^{m}_{1})^{2},\ldots,(\sigma^{m}_{K})^{2})$, let $(\sigma^{m}_{k})^{2} = (t^{m}_{k})^{2}(\sigma^{m}_{dk})^{2}+(1-t^{m}_{k})^{2}(\sigma^{m}_{sk})^{2}$, $\forall k \in [K]$, where demand observations from satellite and drone are  $\hat{\xi}^{m}_{sk} \sim \mathcal{N}(\mu^{m}_{sk},(\sigma^{m}_{sk})^{2})$ and $\hat{\xi}^{m}_{dk} \sim \mathcal{N}(\mu^{m}_{dk},(\sigma^{m}_{dk})^{2})$, respectively. 
        Fig.~\ref{fig:parametric-fusion} shows how we build the parametric data-fusion trust-aided ambiguity set.
        
        We can further extend this process to fuse $H$ data sources. We first assume that $\bm{t} = [\bm{t}_{1},\ldots,\bm{t}_{H}], \bm{t} \in \mathbb{R}^{H\times K}$ with $\sum_{h=1}^{H}t_{hk} = 1$ for all
        $\bm{t}_{k}$, $k \in [K]$. The empirical distribution $\hat{\mathbb{P}}_{N}$ consists of $N$ i.i.d. samples from a Normal distribution $\mathbb{P}^{m}_{e}(\bm{\xi}) = \mathcal{N}(\bm{\mu}^{m},\Sigma^{m})$ after each event $m$. For $\bm{\mu}^{m} = [\mu^{m}_{1},\ldots,\mu^{m}_{K}]^\mathsf{T}$, $\mu^{m}_{k} = \sum_{h=1}^{H} t^{(m-1)}_{hk}\mu^{m}_{hk}$ for all $k \in [K]$; and for $\Sigma^{m} = \text{diag}((\sigma^{m}_{1})^{2},\ldots,(\sigma^{m}_{K})^{2})$, $(\sigma^{m}_{k})^{2} = \sum_{h=1}^{H} (t^{m-1}_{hk})^{2}(\sigma^{m}_{dk})^{2}$ for all $k \in [K]$.

    \subsection{Tractable Reformulation of the DRO Model \eqref{eq:DRO}}\label{2.4}
        The inner worst-case expectation can be rewritten as:
            \begin{equation}\label{eq:DRO-worstcase}
              \sup_{\mathbb{P} \in \mathcal{P}}\mathop{\mathbb{E}}\limits_{\bm{\xi}\sim\mathbb{P}}[\ell(\bm{\xi})], 
            \end{equation}
            where $\ell(\bm{\xi}) = \ell(\bm{x},\bm{\xi})= \sum_{k=1}^{K}\max_{j \in [J]}\ell_{jk}(\xi_{k})$ is a decision-independent separable loss function, and $J$ denotes the number of elementary measurable functions of $\ell(\bm{\xi})$. 
        \begin{theorem}\label{thm:theorem1} 
            If the uncertainty set $\Xi$ is convex and closed, and the loss function is additively separable with respect to $\xi$ and $\{-\ell_{jk}\}_{j \in J}$ is proper, convex, and lower semi-continuous for all $k \in [K]$ \cite{mohajerin2018data}, then \eqref{eq:DRO} is equivalent to \eqref{eq:DRO-reformulation-1}.
        \end{theorem}
    \begin{subequations}\label{eq:DRO-reformulation-1}
            \begin{gather}
                \inf_{\bm{x},\lambda, s_{ik}, \gamma_{ijk}} \lambda\epsilon + \frac{1}{N}\sum^{N}_{i=1}\sum^{K}_{k=1}s_{ik}\label{eq:DRO-reformulation-obj}\\
                \text{s.t.} \quad x \in \mathbb{X},\label{eq:DRO-reformulation-constr1}\\
                b_{jk} + \langle a_{jk},\hat{\xi}_{ik}\rangle + \langle \gamma_{ijk},d_{k} - C_{k}\hat{\xi}_{ik}\rangle \leq s_{ik}, \nonumber\\ \quad\quad\quad i \in [N], j \in [J], k \in [K],\label{eq:DRO-reformulation-constr2}\\
                \|C_{k}^{T}\gamma_{ijk}-a_{jk}\|_{*} \leq \lambda, \forall i \in [N], j \in [J], k \in [K],\label{eq:DRO-reformulation-constr3}\\
                \gamma_{ijk} \geq 0, \quad \forall i \in [N], j \in [J], k \in [K].\label{eq:DRO-reformulation-constr4}
            \end{gather}
    \end{subequations}

\begin{proof}
       Following \eqref{eq:DRO-WassersteinMetric}, we rewrite the worst-case expectation in the DRO model \eqref{eq:DRO} as
    \begin{align}
        \sup_{\mathbb{P} \in \mathcal{P}}\mathop{\mathbb{E}}\limits_{\bm{\xi}\sim\mathbb{P}}[\ell(\bm{\xi})] &=
        \begin{cases}
        \sup_{\Pi, \mathbb{P}} \int_{\Xi} \ell(\bm{\xi}) \mathbb{P}(d\bm{\xi})\\
        \textrm{s.t.} \quad \int_{\Xi^{2}} \|\bm{\xi} - \hat{\bm{\xi}}\| \Pi(d\bm{\xi},d\hat{\bm{\xi}}) \leq \epsilon
        \end{cases}\label{eq:first-equality-1}
    \\&=
        \begin{cases}
        \sup_{\mathbb{P}_{i} \in \mathcal{M}(\Xi)} \frac{1}{N} \sum_{i=1}^{N} \int_{\Xi} \ell(\bm{\xi}) \mathbb{P}_{i}(d\bm{\xi}) \\
        \textrm{s.t.} \quad \frac{1}{N} \sum_{i=1}^{N} \int_{\Xi} \|\bm{\xi} - \hat{\bm{\xi}}_{i}\| \mathbb{P}_{i}(d\bm{\xi}) \leq \epsilon,           
        \end{cases}\label{eq:second-equality-1}
    \end{align}
    where $\Pi$ is the joint distribution of $\bm{\xi}$ and $\hat{\bm{\xi}}$ with marginals $\mathbb{P}$ and $\hat{\mathbb{P}}_{N}$. We drop the minimization problem in the constraint of \eqref{eq:first-equality-1} since the minimization of the Wasserstein metric $d_{W}(\mathbb{P},\hat{\mathbb{P}}_{N})$ is less than equal to radius $\epsilon$ is equivalent as \eqref{eq:first-equality-1} has feasible solution. The second equality \eqref{eq:second-equality-1} indicates that any probability distribution $\Pi$ of $\bm{\xi}$ and $\hat{\bm{\xi}}$ can be constructed from the marginal distribution $\hat{\mathbb{P}}_{N}$ of $\hat{\bm{\xi}}$ and the conditional distribution $\mathbb{P}_{i}$ of $\bm{\xi}$ given $\hat{\bm{\xi}} = \hat{\bm{\xi}}_{i}$, for all $i \in [N]$. Following the standard duality argument \cite{bertsimas1997introduction}, we obtain
    \begin{align}
        & \sup_{\mathbb{P} \in \mathcal{P}}\mathop{\mathbb{E}}\limits_{\xi\sim\mathbb{P}}[\ell(\bm{\xi})] =
        \sup_{\mathbb{P}_{i} \in \mathcal{M}(\Xi)} \quad \Big\{\inf_{\lambda \geq 0} \frac{1}{N} \int_{\Xi} \ell(\bm{\xi}) \mathbb{P}_{i}(d\bm{\xi}) \\ \nonumber
        & + \lambda (\epsilon - \frac{1}{N}\sum_{i=1}^{N}\int_{\Xi} \|\bm{\xi} - \hat{\bm{\xi}}_{i}\| \mathbb{P}_{i}(d\bm{\xi}))\Big\} \\
        & \leq
        \inf_{\lambda \geq 0} \quad \Big \{ \sup_{\mathbb{P}_{i} \in \mathcal{M}(\Xi)} \lambda\epsilon + \frac{1}{N}\sum_{i=1}^{N}\int_{\Xi}(\ell(\bm{\xi}) - \lambda\|\bm{\xi} - \hat{\bm{\xi}}_{i}\| \mathbb{P}_{i}(d\bm{\xi}))\Big\} \label{eq:max-min-1}\\
        & =
        \inf_{\lambda \geq 0} \quad \lambda\epsilon + \frac{1}{N}\sum_{i=1}^{N} \sup_{\bm{\xi} \in \Xi}  \big\{ \ell(\bm{\xi}) - \lambda\|\bm{\xi} - \hat{\bm{\xi}}_{i}\| \big\},\label{eq:Dirac-1}
    \end{align}
    where \eqref{eq:max-min-1} holds because of the max-min inequality, and \eqref{eq:Dirac-1} follows from the fact that $\mathcal{M}(\Xi)$ contains all the Dirac distributions supported on $\Xi$.
    Meanwhile, the loss function in our problem is additively separable with respect to the temporal structure of $\bm{\xi}$, that is,
    \begin{equation*}
        \ell(\bm{\xi}):=\sum^{K}_{k=1}\max_{j \in [J]} \ell_{jk}(\xi_{k}),
    \end{equation*}
    where $\ell_{jk} : \mathbb{R} \rightarrow \bar{\mathbb{R}}$ is a measurable function for any $j \in [J]$ and $k \in [K]$. Since we use $L_1$-norm to define the Wasserstein metric, $\|\cdot\|_{\text{K}}$ reduces to $L_1$-norm on $\mathbb{R}^{K}$. Now, \eqref{eq:Dirac-1} can be written with the interchange of the summation and the maximization, which yields 
    \begin{align}
        \sup_{\mathbb{P} \in \mathcal{P}}\mathop{\mathbb{E}}\limits_{\bm{\xi}\sim\mathbb{P}}[\ell(\bm{\xi})] &=
        \inf_{\lambda \geq 0} \quad \lambda\epsilon + \frac{1}{N}\sum_{i=1}^{N} \sup_{\bm{\xi} \in \Xi} (\ell(\bm{\xi}) - \lambda\|\bm{\xi} - \hat{\bm{\xi}}_{i}\|)\nonumber \\
        &=
        \inf_{\lambda \geq 0} \quad \lambda\epsilon + \frac{1}{N} \sum_{i-1}^{N} \sum_{k=1}^{K} \sup_{\xi_{k} \in \Xi_{k}} (\max_{j=1,\ldots,J} \ell_{jk}(\xi_{k}) \\ \nonumber
            & - \lambda \|\xi_{k} - \hat{\xi}_{ik}\|)\label{eq:interchange-1}.
    \end{align}
    After introducing auxiliary variables in \eqref{eq:interchange-1}, we have
    \begin{align}
        & \quad
        \begin{cases}
            \inf_{\lambda,s_{ik}} \quad \lambda\epsilon + \frac{1}{N}\sum_{i=1}^{N}\sum_{k=1}^{K}s_{ik}\\
            \textrm{s.t.} \quad \sup_{\xi_{k} \in \Xi_{k}}(\ell_{jk}(\xi_{k})-\lambda\|\xi_{k} - \hat{\xi}_{ik}\|) \leq s_{ik} \\ \quad\quad\quad \forall i \in [N], j \in [J], k \in [K]\\
            \quad \quad \lambda \geq 0
        \end{cases}\\
        & \leq
        \begin{cases}
            \inf_{\lambda,s_{ik},z_{ijk}} \quad \lambda\epsilon + \frac{1}{N}\sum_{i=1}^{N}\sum_{k=1}^{K}s_{ik}\\
            \textrm{s.t.} \quad \sup_{\xi_{k} \in \Xi_{k}}(\ell_{jk}(\xi_{k})-\langle z_{ijk},\xi_{k}\rangle)+\langle z_{ijk},\hat{\xi}_{ik}\rangle \leq s_{ik} \\ \quad\quad\quad \forall i \in [N], j \in [J], k \in [K]\\
            \quad \quad \|z_{ijk}\|_{*} \leq \lambda \quad \forall i \in [N], j \in [J], k \in [K]\\
        \end{cases}\\
        & =
        \begin{cases}
            \inf_{\lambda,s_{ik},z_{ijk}} \quad \lambda\epsilon + \frac{1}{N}\sum_{i=1}^{N}\sum_{k=1}^{K}s_{ik}\\
            \textrm{s.t.} \quad [-\ell_{jk} + \chi_{\Xi_{k}}]^{*}(-z_{ijk})+\langle z_{ijk},\hat{\xi}_{ik}\rangle \leq s_{ik} \\ \quad\quad\quad \forall i \in [N], j \in [J], k \in [K]\\
            \quad \quad \|z_{ijk}\|_{*} \leq \lambda \quad \forall i \in [N], j \in [J], k \in [K],\\
        \end{cases}        
    \end{align}
    where the inequality holds as an equality provided that $\Xi_{k}$ and $\{\ell_{jk}\}_{j \in [J]}$ satisfy the convexity assumption for all $k \in [K]$. Finally, by \cite{rockafellar2009variational}, the conjugate of $-\ell_{jk} + \chi_{\Xi_{k}}$ can be replaced by the inf-convolution of the conjugates of $-\ell_{jk}$ and $\chi_{\Xi_{k}}$. By definition of the conjugacy operator, we have
    \begin{align}
        [-\ell_{jk}]^{*}(z) &= [-a_{jk}]^{*}(z) = \sup_{\xi}\langle z,\xi_{k} \rangle + \langle a_{jk},\xi_{k} \rangle + b_{jk}\nonumber\\
        &= \begin{cases}
    b_{jk} & \text{if}\quad z = -a_{jk},\\\nonumber
    \infty & \text{else},
    \end{cases}
    \end{align}
    and
    \[\sigma_{\Xi_{k}(\nu)} =
    \begin{cases}
    \sup_{\xi_{k}} \quad \langle \nu,\xi_{k} \rangle \\
    \text{s.t.} \quad C_{k}\xi_{k} \leq d_{k}
    \end{cases} = 
    \begin{cases}
    \inf_{\gamma \geq 0} \quad \langle \gamma,d \rangle\\
    \text{s.t.} \quad C_{k}^{T}\gamma = \nu,
    \end{cases}\]
where the last equality follows from strong duality, which holds as the uncertainty set is non-empty.
After bringing this form to \eqref{eq:DRO}, we obtain the equivalent linear programming reformulation \eqref{eq:DRO-reformulation-1}. 
\end{proof}

\section{Trust Update and Trust Selection}\label{3}

      \subsection{Parametric Data-fusion Trust Update Procedures}\label{3.1}

           The complete trust update process is illustrated as Alg.~\ref{alg:trust-update1}. Before the process begins, we hold an original trust $\bm{t}^{0}_{h} = (t^{0}_{h1},\ldots,t^{0}_{hK})^\mathsf{T}$ over each information source $h$, $\forall h \in [H]$. For each data event $m$, we will predict a distribution $\mathbb{P}^{m}_{h}$ for each source $h$. After event $m$ ends, we then observe the realization $\bm{\xi}^{m}_{true}$ of the uncertain demand. We assume that the relative error $\bm{r}_{h}$ between the mean value of the predicted distribution provided by source $h$ and the mean value of true demand distribution is fixed, which means $\bm{\mu}_{h} = \langle \bm{\mu}_{true},\bm{r}_{h} \rangle$. For each event $m$, when $ \mathbb{P}_{true}$ changes, $\mathbb{P}_{h}$ will change based on the relative error $\bm{r}_{h}$. Therefore, as the number of events (i.e., $M$) grows, we should see a trend of the trust value $t^{m}_{hk}$ in source $h$ at region $k$ getting closer to and fluctuating near $t^{ideal}_{hk}$, which satisfies $\mu_{true,k} = \sum_{h=1}^{H}(t^{ideal}_{hk}\times \mu_{hk})$, but might not exactly equals $t^{ideal}_{hk}$. Once the trust value starts to fluctuate within a certain range, we assume that it has reached the fluctuating interval.
   
\begin{figure}[htbp]	
            \centering
            \includegraphics[width=0.4\textwidth]{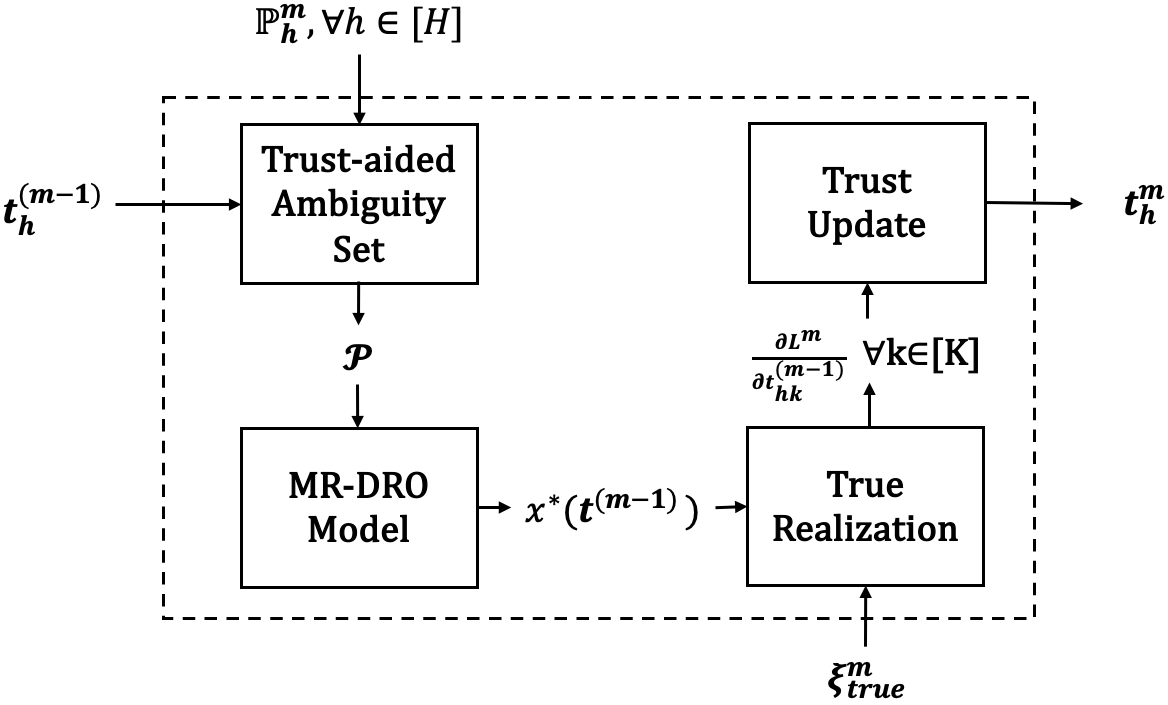}
            \caption{Illustrating parametric data-fusion trust update based on losses}	
            \label{fig:parafusion-trust-update-step}
        \end{figure}

        We view each event as a 4-step procedure, which is shown in Fig.~\ref{fig:parafusion-trust-update-step}.    

At the beginning of a wildfire event $m$, we hold trust $\bm{t}^{(m-1)}_{h}$ in source $h$ and receive the predicted distribution $\mathbb{P}^{m}_{h}$. Based on the current trust and data, we can construct the trust-aided ambiguity set $\mathcal{P}$. Then we solve the MR-DRO model to attain an optimal resource allocation solution $\bm{x}^{*}(\bm{t}^{(m-1)})$ based on the current trust. 
        
        We update the trust based on losses calculated as follows. After event $m$ ends, we obtain information about the true demand distribution $\bm{\xi}^{m}_{true}$. We calculate the corresponding loss $L^{m}$ through:       
        \begin{align}\label{eq:real}
        L^{m} & = \ell(\bm{x}^{*},\bm{\xi}^{m}_{true})\\
        & = (\bm{c}^{u})^{\mathsf{T}}[\bm{\xi}^{m}_{true}-\bm{x}^{*}(\bm{t}^{(m-1)})]^{+}+(\bm{c}^{o})^{\mathsf{T}}[\bm{x}^{*}(\bm{t}^{(m-1)})-\bm{\xi}^{m}_{true}]^{+}\nonumber\\
        & = \sum_{k=1}^{K}[c^{u}_{k}[\xi^{m}_{true,k}-x_{k}^{*}(\bm{t}^{(m-1)})]^{+}+c^{o}_{k}[x_{k}^{*}(\bm{t}^{(m-1)})-\xi^{m}_{true,k}]^{+}]. \nonumber
        \end{align}
We then update each $t^{m}_{hk}$ for all $k \in [K]$ in $\bm{t}^{m}_{h}$ based on the partial derivative $\frac{\partial L^{m}}{\partial t^{(m-1)}_{k}}$:
        \begin{equation}\label{eq:trust1}
                t^{m}_{k} = t^{(m-1)}_{k} - w \times \frac{\partial L^{m}}{\partial t^{(m-1)}_{k}},
        \end{equation}
        where $w$ is a small step size; the negative partial derivative of losses in event $m$ with respect to trust $t^{(m-1)}_{k}$ at region $k$ indicates that an increment of the trust $t^{(m-1)}_{k}$ will result in a decrease of real loss, for which we increase the trust $t^{m}_{k}$ for the next event. 
            
    \section{Numerical Studies}\label{4}
    We conduct numerical tests and compare different models and approaches using a diverse set of wildfire suppression instances. We use Gurobi 9.5.2 for solving all linear programming models. The algorithm for trust update is implemented in Python 3.9.12. All numerical tests are conducted on a MacBook Pro with 16 GB RAM and an Apple M1 Pro chip.
    
    \subsection{Experimental Design}\label{4.1}
        We first consider a baseline case with $K=3$, $H=2$, and $M=50$. We use $N=200$ samples and set $\epsilon = 0.01$ as the radius of the ambiguity set $\mathcal{P}$. The SP model \eqref{eq:SP} uses all $|\Omega|= N =200$ samples. The unit penalty costs for over-served and unmet demand are $\bm{c}^{u} = (5000,5000,5000)^\mathsf{T}$ and $\bm{c}^{o} = (1000,1000,1000)^\mathsf{T}$, respectively. The resource budget is $B = 1000$. The true demand for region $k$ is an integer uniformly sampled from $[100,200]$, for all $k \in [K]$. We let $\bm{r}_{h_{1}} = (1.1,0.6,1.1)^\mathsf{T}$, $\bm{r}_{h_{2}} = (0.7,1.2,0.3)^\mathsf{T}$, $\bm{\sigma}_{s} = 0.02*\bm{\mu}_{s}$, and $\bm{\sigma}_{d} = 0.02*\bm{\mu}_{d}$. We set the original trust of drone and satellite data being the same, such that $\bm{t}_{0} = (0.5,0.5,0.5)^\mathsf{T}$. In the out-of-sample tests, we consider $Q = 100$ events. Due to the complexity of the loss function \eqref{eq:real}, in numerical experiments, we use the incremental trial-and-error method to judge the partial derivative with respect to trust with $w= 10^{-3}$.

                     \begin{algorithm}[htbp]
\caption{Parametric Data-fusion trust update}\label{alg:trust-update1}
    \begin{algorithmic}[1]
    \State Inputs: original trust $\bm{t}^{0} = (\bm{t}^{0}_{1},\ldots,\bm{t}^{0}_{H})$, step size $w$.
    \For{$m=1,\ldots, M$}
        \State Generate the trust-aided ambiguity set $\mathcal{P}$ with $\bm{t}^{(m-1)}$.
        \State Solve reformulation \eqref{eq:DRO-reformulation-1} of the MR-DRO model with $\mathcal{P}$ and obtain optimal solution $\bm{x}^{*}(\bm{t}^{(m-1)})$.
        \For{$k=1,\ldots,K$}
            \For {$h=1,\ldots,H$}
            \State Generate $t^{(m-1)}_{h^{increase}k}$ and $t^{(m-1)}_{h^{decrease}k}$.
            \State Solve \eqref{eq:DRO-reformulation-1} with ambiguity set based on $t^{(m-1)}_{h^{increase}k}$ and $t^{(m-1)}_{h^{decrease}k}$.
            \State Calculate loss $L_{h^{increase}}$ and $L_{h^{decrease}}$.
            \State Estimate partial derivative $\frac{\partial L}{\partial t^{(m-1)}_{k}}$ at $t^{(m-1)}_{hk}$.
            \State Let $(t^{new})^{(m-1)}_{hk}=t^{(m-1)}_{hk} - \frac{\partial L}{\partial t^{(m-1)}_{k}}\times w$.
            \EndFor
            \For {$h=1,\ldots,H$}
                \State Normalize $(t^{new})^{(m-1)}_{hk}$.
            \EndFor
        \EndFor
    \State Set $\bm{t}^{m} \leftarrow (\bm{t}^{new})^{(m-1)}$.
    \EndFor
\end{algorithmic}
\end{algorithm}

    \subsection{Trust Update and Computational Results}\label{4.2}
\begin{figure}[htbp]
        \centering
        \includegraphics[width=0.4\textwidth]{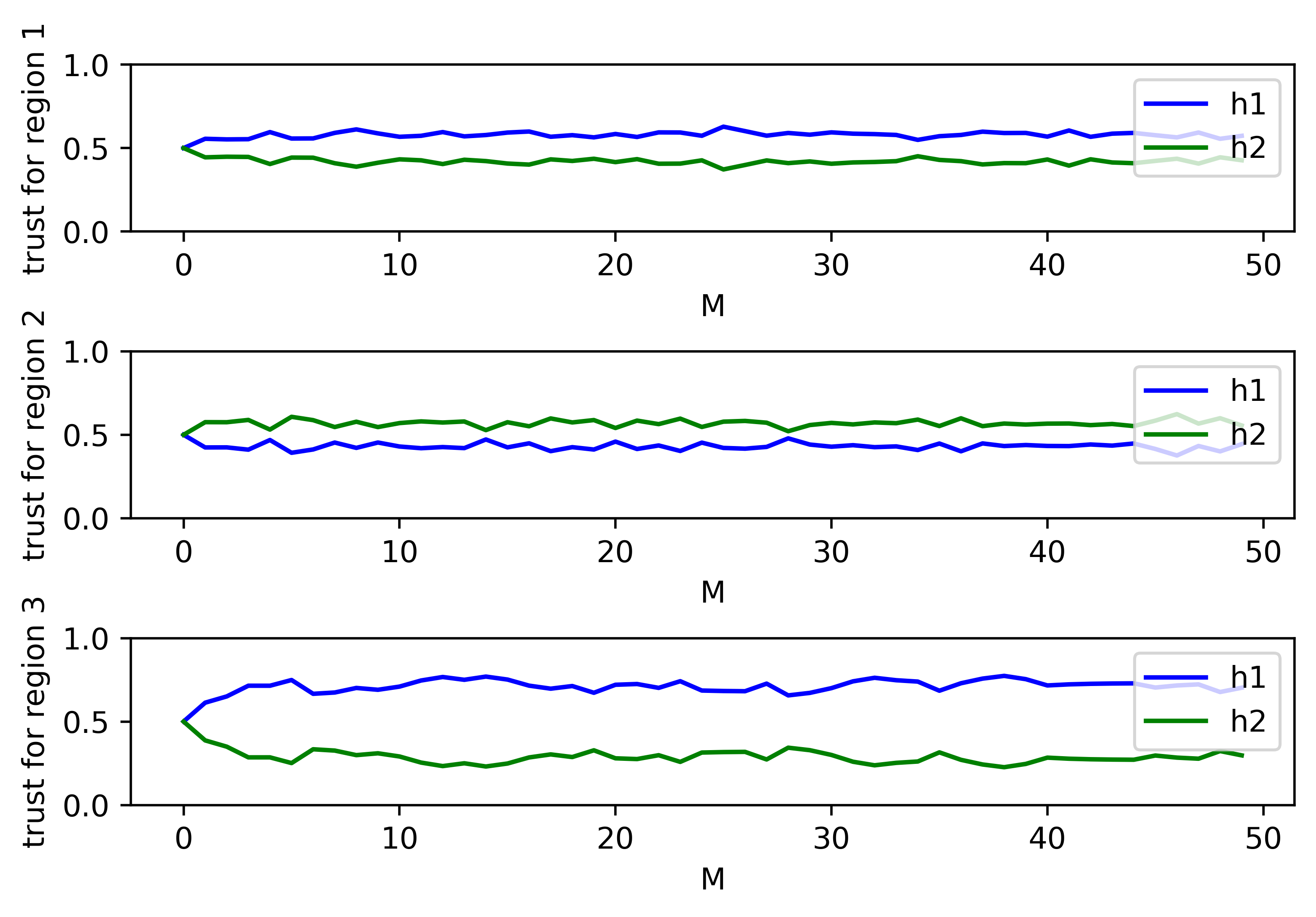}
        \caption{Trust update process with the baseline setting}
        \label{fig:baseline}
    \end{figure}

    Using the baseline setting and the trust update algorithm, we first obtain the result of the trust update process, reported in Table \ref{tab:trust-update-result} and Fig.~\ref{fig:baseline}. In Table \ref{tab:trust-update-result}, ``Trust Interval,'' ``Loss'' and ``Time'' denote the range of fluctuating trust, the average losses with trust in the fluctuating interval in thousand dollars and the total computation time in seconds, respectively.
    
    \begin{table}[htbp]
    \begin{center}
    \caption{Trust update in the baseline setting}
    \resizebox{0.48\textwidth}{!}{
    \begin{tabular}{cccccc} 
     \toprule
     $K$ & $M$ & Source & Trust Interval & Loss (*\$1000) & Time (sec.)\\
     \midrule
     \multirow{ 6}{*}{$3$} & \multirow{ 6}{*}{$50$} & \multirow{ 3}{*}{$h_{1}$ = drone} &$[0.55,0.63]$ & \multirow{ 6}{*}{$152.25$} & \multirow{ 6}{*}{$155.40$}\\
     &  &  & $[0.38,0.48]$ &  & \\
     &  &  & $[0.66,0.77]$ &  & \\
     &  & \multirow{ 3}{*}{$h_{2}$ = satellite} & $[0.37,0.45]$ &  & \\
     &  &  & $[0.52,0.62]$ &  & \\
     &  &  & $[0.23,0.34]$ &  & \\
    \bottomrule
    \end{tabular}
    }
    \label{tab:trust-update-result}
    \end{center}
    \end{table}

\begin{figure}[htbp]
        \centering
        \includegraphics[width=0.4\textwidth]{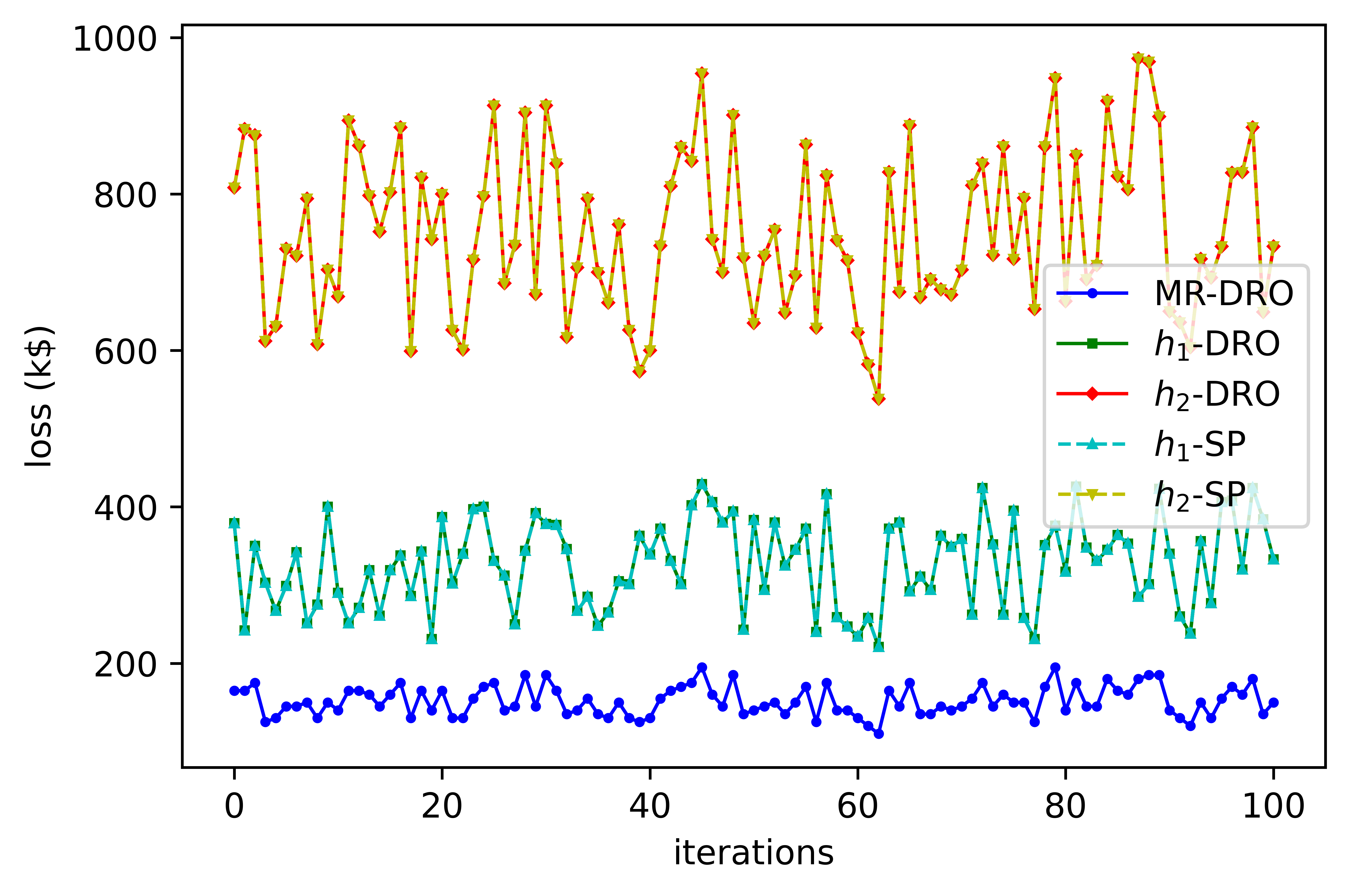}
        \caption{Out-of-sample performances with different approaches ($B = 1000$)}
        \label{fig:out-of-sample-comparison1}
    \end{figure}

    We vary the trust in different regions to consider two trust vectors $\bm{t}^{*}_{h_{1}} = (0.58,0.43,0.72)^\mathsf{T}$ and $\bm{t}^{*}_{h_{2}} = (0.42,0.57,0.28)^\mathsf{T}$. We then examine out-of-sample performances of MR-DRO model with two information sources, a DRO and an SP model based on one single information source (either drone or satellite), and report the results in Table \ref{tab:comparison1} and Fig.~\ref{fig:out-of-sample-comparison1}. The true demand distribution varies in each event of the out-of-sample tests, with the mean value of the true marginal demand distribution for region $k$ in event $h$ as a random integer between [100,200], for all $k \in [K]$ and $q \in [Q]$. The relative error relationship between the predicted distribution of the satellite or the drone and the true distribution keeps the same with the baseline setting.

    In Table \ref{tab:comparison1}, the DRO and SP with source as prefix denotes DRO or SP models solved with information from a single source. Note that MR-DRO model has the best performance among all DRO models in terms of the average loss, with similar computational time. The reason why MR-DRO model performs better than $h_{1}$-DRO and $h_{2}$-DRO is that based on the trust we obtain from previous trust update process, we are able to reduce the effect of prediction error caused by $h_{1}$ and $h_{2}$. In addition, given sufficient budget and with current radius $\epsilon$, the solution offered by single source DRO and single source SP are the same.
    \begin{table}[htbp]
    \begin{center}
    \caption{Out-of-sample performances of MR-DRO and single-sourced DRO and SP with budget $B = 1000$}
    \begin{tabular}{cccccc} 
     \toprule
     $K$ & $M$ & $Q$ & Method & Loss (*\$1000) & Time (sec.)\\
     \midrule
     \multirow{ 5}{*}{3} & \multirow{ 5}{*}{50} & \multirow{ 5}{*}{100}& MR-DRO & 155.79 & 18.49\\
      &  &  & $h_1$-DRO & 325.11 & 18.35\\
      &  &  & $h_2$-DRO & 775.03 & 18.98\\
      &  &  & $h_1$-SP & 325.11 & 5.48\\
      &  &  & $h_2$-SP & 775.03 & 5.48\\
    \bottomrule
    \end{tabular}
    \label{tab:comparison1}
    \end{center}
    \end{table}
    
    \subsection{Sensitivity Analysis}\label{4.3}
    \subsubsection{Varying Budget $B$}
\begin{figure}[htbp]
        \centering
        \includegraphics[width=0.4\textwidth]{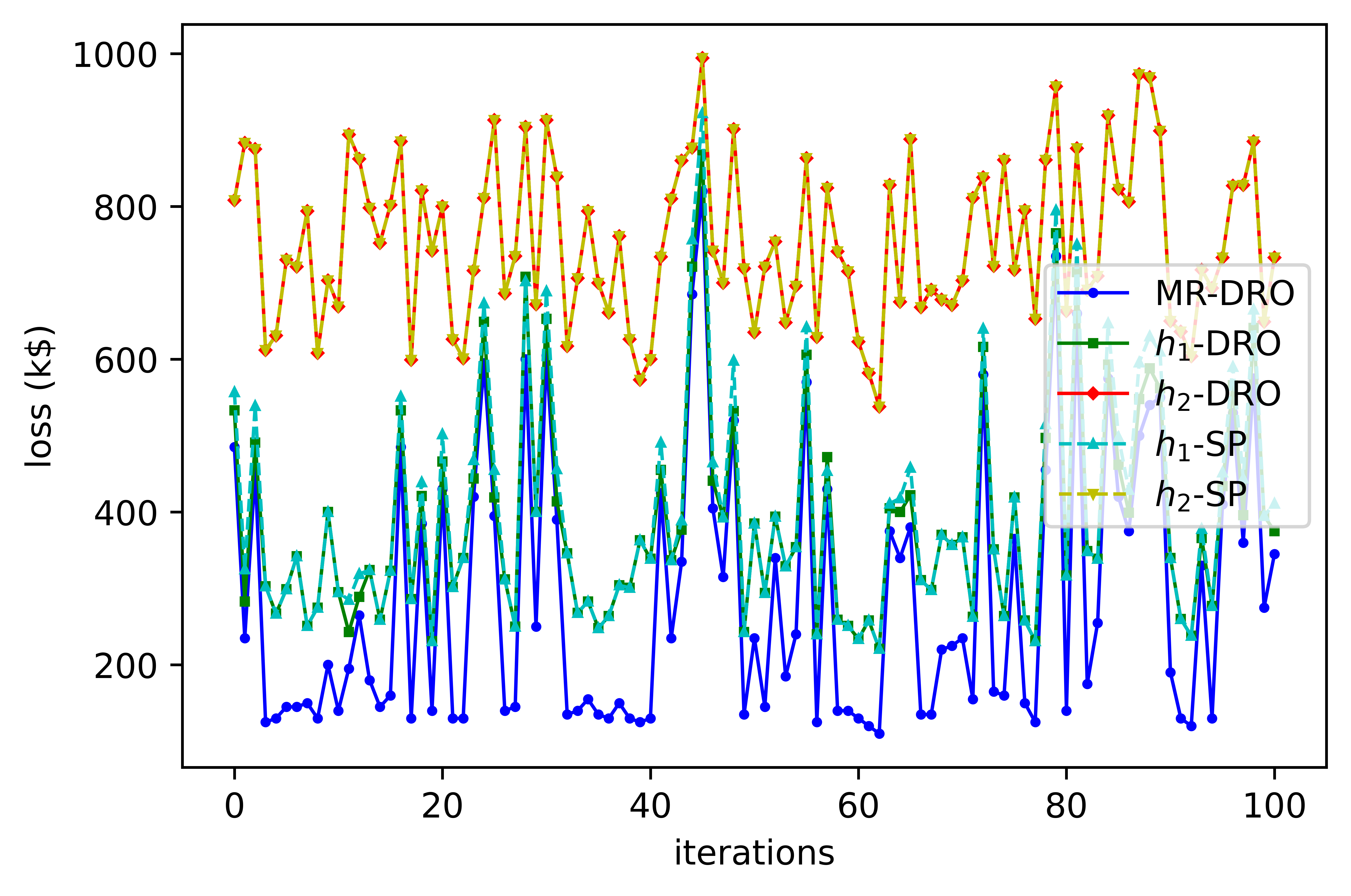}
        \caption{Out-of-sample performances with different approaches ($B = 400$)}
        \label{fig:out-of-sample-comparison2}
    \end{figure}

    We keep other parameters same as the baseline setting and change the budget to $B = 400$ to investigate how MR-DRO model and other methods perform when the resource allocation budget becomes tight. The results are reported in Table \ref{tab:comparison2} and Fig.~\ref{fig:out-of-sample-comparison2}, showing that whether the budget is always sufficient or occasionally insufficient, MR-DRO model generally performs better than DRO and SP method with single source in terms of average loss.
    \begin{table}[htbp]
    \begin{center}
    \caption{Out-of-sample performances of MR-DRO and single-sourced DRO and SP with budget $B = 400$}
    \begin{tabular}{cccccc} 
     \toprule
     $K$ & $M$ & $Q$ & Method & Loss (*\$1000) & Time (sec.)\\
     \midrule
     \multirow{ 5}{*}{3} & \multirow{ 5}{*}{50} & \multirow{ 5}{*}{100}& MR-DRO & 304.60 & 17.33\\
      &  &  & $h_1$-DRO & 394.61 & 17.47\\
      &  &  & $h_2$-DRO & 776.74 & 18.55\\
      &  &  & $h_1$-SP & 409.88 & 5.48\\
      &  &  & $h_2$-SP & 777.10 & 5.60\\
    \bottomrule
    \end{tabular}
    \label{tab:comparison2}
    \end{center}
    \end{table}
    
    \subsubsection{Varying Sample Size $M$}

    We conduct numerical experiments to see if the variation of $M$ would affect the trust value we get, reported in Table \ref{tab:trust-update-varying-M}. When the number of events in the trust update process is extremely small (i.e., $M=10$), the trust values cannot reach a fluctuating interval. As the number of events increases, the $t^{*}$-value for each source does not vary significantly, which indicates that we are able to attain a stabilized trust value in finite steps.

    \subsection{Computational Scalability Analysis}\label{4.4}
    
    In this section, we vary $K$ from 3 to 5 and 10 to see how the number of regions will affect the computational time. Table \ref{tab:trust-update-varying-K} shows that as the number of regions increases, the computational time increases significantly. This is because the larger number of regions not only leads to an increase in computational time for solving the MR-DRO model with fixed trust, but also results in longer time of the trust update process.

      \begin{table}[htbp]
    \begin{center}
    \caption{Trust update process with varying $M$}
    \begin{tabular}{cccccc} 
     \toprule
     $K$ & $M$ & $t_{h_{1}}$ Interval & $t_{h_{2}}$& $t^{*}_{h_{1}}$ & $t^{*}_{h_{2}}$\\
     \midrule
     \multirow{ 15}{*}{$3$} & \multirow{ 3}{*}{$10$} & N/A & N/A & N/A & N/A\\
      & & N/A & N/A & N/A & N/A\\
      & & N/A & N/A & N/A & N/A\\
    & \multirow{ 3}{*}{$50$} & [0.55, 0.63] & [0.37, 0.45] & 0.58 & 0.42\\
    & & [0.38, 0.48] & [0.52, 0.62] & 0.43 & 0.57\\
    & & [0.66, 0.77] & [0.23, 0.34] & 0.72 & 0.28\\
    & \multirow{ 3}{*}{$100$} & [0.53, 0.63] & [0.37, 0.47] & 0.58 & 0.42\\
    & & [0.38, 0.48] & [0.52, 0.62] & 0.43 & 0.57\\
    & & [0.66, 0.79] & [0.21, 0.34] & 0.72 & 0.28\\
    & \multirow{ 3}{*}{$150$} & [0.53, 0.64] & [0.36, 0.47] & 0.58 & 0.42\\
    & & [0.36, 0.49] & [0.51, 0.64] & 0.43 & 0.57\\
    & & [0.66, 0.83] & [0.17, 0.34] & 0.72 & 0.28\\
    & \multirow{ 3}{*}{$200$} & [0.52, 0.64] & [0.36, 0.48] & 0.58 & 0.42\\
    & & [0.36, 0.49] & [0.51, 0.64] & 0.43 & 0.57\\
    & & [0.66, 0.83] & [0.17, 0.34] & 0.72 & 0.28\\
    \bottomrule
    \end{tabular}
    \label{tab:trust-update-varying-M}
    \end{center}
    \end{table}

    \begin{table}[htbp]
    \begin{center}
    \caption{Trust Update process with varying $K$}
    \begin{tabular}{ccc} 
     \toprule
     $K$ & $M$ & Time (sec.)\\
     \midrule
     3 & \multirow{ 3}{*}{200} & 623.67\\
     5 &  &1571.55\\
     10 &  &6621.60\\
    \bottomrule
    \end{tabular}
    \label{tab:trust-update-varying-K}
    \end{center}
    \end{table}

    \section{Conclusions}\label{5}
    In this work, we formulated a MR-DRO model for solving the stochastic resource allocation problem, where we constructed a parametric data-fusion trust-aided ambiguity set for the DRO model. We also proposed a trust update process, to dynamically adjust trust after realizing the uncertainty and collecting more data from different sources. Our numerical results indicated that the MR-DRO model performs better than solving DRO or SP model by using only one information source.

    The results have the following limitations. First, the current way of information fusion in our model for constructing the trust-aided ambiguity set requires both the predicted distribution and the true distribution to follow Normal distributions, which might not be true in specific applications. Second, the proposed trust update process assumes that the prediction ability of different sources is constant across all the events. Third, we model the wildfire suppression resource allocation problem as a static one-stage problem. To address the three limitations, future research is needed to find a different way of data-fusion for constructing the ambiguity set, which does not require strong assumptions on the type of distributions. Similarly, one can verify whether the current way of trust update can work for situations where the prediction ability of the satellite and the drone are not fixed. Another promising future research direction is to model the wildfire suppression resource allocation problem as a multi-stage stochastic dynamic program, which may introduce additional computational challenges.


\begin{thebibliography}{10}
\providecommand{\url}[1]{#1}
\csname url@samestyle\endcsname
\providecommand{\newblock}{\relax}
\providecommand{\bibinfo}[2]{#2}
\providecommand{\BIBentrySTDinterwordspacing}{\spaceskip=0pt\relax}
\providecommand{\BIBentryALTinterwordstretchfactor}{4}
\providecommand{\BIBentryALTinterwordspacing}{\spaceskip=\fontdimen2\font plus
\BIBentryALTinterwordstretchfactor\fontdimen3\font minus
  \fontdimen4\font\relax}
\providecommand{\BIBforeignlanguage}[2]{{%
\expandafter\ifx\csname l@#1\endcsname\relax
\typeout{** WARNING: IEEEtran.bst: No hyphenation pattern has been}%
\typeout{** loaded for the language `#1'. Using the pattern for}%
\typeout{** the default language instead.}%
\else
\language=\csname l@#1\endcsname
\fi
#2}}
\providecommand{\BIBdecl}{\relax}
\BIBdecl

\bibitem{delage2010distributionally}
E.~Delage and Y.~Ye, ``Distributionally robust optimization under moment
  uncertainty with application to data-driven problems,'' \emph{Operations
  Research}, vol.~58, no.~3, pp. 595--612, 2010.

\bibitem{jiang2016data}
R.~Jiang and Y.~Guan, ``Data-driven chance constrained stochastic program,''
  \emph{Mathematical Programming}, vol. 158, no.~1, pp. 291--327, 2016.

\bibitem{mohajerin2018data}
P.~Mohajerin~Esfahani and D.~Kuhn, ``Data-driven distributionally robust
  optimization using the wasserstein metric: Performance guarantees and
  tractable reformulations,'' \emph{Mathematical Programming}, vol. 171, no.
  1-2, pp. 115--166, 2018.

\bibitem{yang2021multi}
M.~Yang, Y.~Liu, and G.~Yang, ``Multi-period dynamic distributionally robust
  pre-positioning of emergency supplies under demand uncertainty,''
  \emph{Applied Mathematical Modelling}, vol.~89, pp. 1433--1458, 2021.

\bibitem{wang2021two}
W.~Wang, K.~Yang, L.~Yang, and Z.~Gao, ``Two-stage distributionally robust
  programming based on worst-case mean-cvar criterion and application to
  disaster relief management,'' \emph{Transportation Research Part E: Logistics
  and Transportation Review}, vol. 149, p. 102332, 2021.

\bibitem{basciftci2023resource}
B.~Basciftci, X.~Yu, and S.~Shen, ``Resource distribution under spatiotemporal
  uncertainty of disease spread: Stochastic versus robust approaches,''
  \emph{Computers \& Operations Research}, vol. 149, p. 106028, 2023.

\bibitem{durrant2016multisensor}
H.~Durrant-Whyte and T.~C. Henderson, ``Multisensor data fusion,''
  \emph{Springer handbook of robotics}, pp. 867--896, 2016.

\bibitem{abdulhafiz2014bayesian}
W.~A. Abdulhafiz and A.~Khamis, ``Bayesian approach with pre-and post-filtering
  to handle data uncertainty and inconsistency in mobile robot local
  positioning,'' \emph{Journal of Intelligent Systems}, vol.~23, no.~2, pp.
  133--154, 2014.

\bibitem{dietrich2017probabilistic}
F.~Dietrich and C.~List, ``Probabilistic opinion pooling generalized. part one:
  general agendas,'' \emph{Social Choice and Welfare}, vol.~48, no.~4, pp.
  747--786, 2017.

\bibitem{noyvirt2012human}
A.~Noyvirt and R.~Qiu, ``Human detection and tracking in an assistive living
  service robot through multimodal data fusion,'' in \emph{IEEE 10th
  International Conference on Industrial Informatics}.\hskip 1em plus 0.5em
  minus 0.4em\relax IEEE, 2012, pp. 1176--1181.

\bibitem{howcroft2017prospective}
J.~Howcroft, J.~Kofman, and E.~D. Lemaire, ``Prospective fall-risk prediction
  models for older adults based on wearable sensors,'' \emph{IEEE Transactions
  on Neural Systems and Rehabilitation Engineering}, vol.~25, no.~10, pp.
  1812--1820, 2017.

\bibitem{pan2020human}
D.~Pan, H.~Liu, D.~Qu, and Z.~Zhang, ``Human falling detection algorithm based
  on multisensor data fusion with svm,'' \emph{Mobile Information Systems},
  vol. 2020, pp. 1--9, 2020.

\bibitem{bloomfield2021machine}
R.~A. Bloomfield, J.~S. Broberg, H.~A. Williams, B.~A. Lanting, K.~A. McIsaac,
  and M.~G. Teeter, ``Machine learning and wearable sensors at preoperative
  assessments: functional recovery prediction to set realistic expectations for
  knee replacements,'' \emph{Medical Engineering \& Physics}, vol.~89, pp.
  14--21, 2021.

\bibitem{nsengiyumva2022critical}
W.~Nsengiyumva, S.~Zhong, M.~Luo, Q.~Zhang, and J.~Lin, ``Critical insights
  into the state-of-the-art nde data fusion techniques for the inspection of
  structural systems,'' \emph{Structural Control and Health Monitoring},
  vol.~29, no.~1, p. e2857, 2022.

\bibitem{lee2004trust}
J.~D. Lee and K.~A. See, ``Trust in automation: Designing for appropriate
  reliance,'' \emph{Human factors}, vol.~46, no.~1, pp. 50--80, 2004.

\bibitem{guo2021modeling}
Y.~Guo and X.~J. Yang, ``Modeling and predicting trust dynamics in human--robot
  teaming: A bayesian inference approach,'' \emph{International Journal of
  Social Robotics}, vol.~13, no.~8, pp. 1899--1909, 2021.

\bibitem{yang2023toward}
X.~J. Yang, C.~Schemanske, and C.~Searle, ``Toward quantifying trust dynamics:
  How people adjust their trust after moment-to-moment interaction with
  automation,'' \emph{Human Factors}, vol.~65, no.~5, pp. 862--878, 2023.

\bibitem{kleywegt2002sample}
A.~J. Kleywegt, A.~Shapiro, and T.~Homem-de Mello, ``The sample average
  approximation method for stochastic discrete optimization,'' \emph{SIAM
  Journal on Optimization}, vol.~12, no.~2, pp. 479--502, 2002.

\bibitem{bertsimas1997introduction}
D.~Bertsimas and J.~N. Tsitsiklis, \emph{Introduction to linear
  optimization}.\hskip 1em plus 0.5em minus 0.4em\relax Athena Scientific
  Belmont, MA, 1997, vol.~6.

\bibitem{rockafellar2009variational}
R.~T. Rockafellar and R.~J.-B. Wets, \emph{Variational analysis}.\hskip 1em
  plus 0.5em minus 0.4em\relax Springer Science \& Business Media, 2009, vol.
  317.

\end{thebibliography}
   

\end{document}